\documentclass{article}
\usepackage{spconf,amsmath}
\usepackage[]{hyperref}
\usepackage{mathptmx, amsmath, amssymb, amsthm, mathptmx, enumerate, color}
\usepackage{graphicx}
\usepackage{epstopdf}

\usepackage{mathrsfs}

\ninept
\newtheorem{theorem}{Theorem}[section]

\newtheorem{proposition}{Proposition}[section]
\theoremstyle{definition}

\newtheorem{example}{Example}[section]
\newtheorem{remark}{Remark}[section]

\numberwithin{equation}{section}

\newcommand{\ip}[2]{\left\langle #1 , #2 \right\rangle}
\newcommand{\norm}[1]{\left\| #1 \right\|}
\newcommand{\minimize}[0]{\operatornamewithlimits{minimize\ }}

\newcommand{\dom}{\operatorname{dom}}

\newcommand{\lrangle}[1]{\left\langle #1 \right\rangle}
\newcommand{\argmin}{\operatornamewithlimits{argmin\ }}
\newcommand{\Id}{\operatorname{Id}}
\newcommand{\prox}{\operatorname{Prox}}
\newcommand{\fix}{\operatorname{Fix}}
\newcommand{\setLO}[2]{\mathcal{B}\left( #1 , #2 \right)}
\newcommand{\zeroMatrix}{\mathrm{O}}

\newcommand{\setN}{\mathbb{N}}
\newcommand{\setR}{\mathbb{R}}

\newcommand{\transT}{\mathsf{T}}

\newcommand{\spX}{\mathcal{X}}
\newcommand{\spY}{\mathcal{Y}}
\newcommand{\spZ}{\mathcal{Z}}
\newcommand{\sptildeZ}{\widetilde{\mathcal{Z}}}
\newcommand{\spfrakZ}{\mathfrak{Z}}
\newcommand{\spH}{\mathcal{H}}
\newcommand{\spK}{\mathcal{K}}

\newcommand{\opL}{\mathfrak{L}}

\newcommand{\opC}{\mathfrak{C}}

\newcommand{\xstar}{{x^\star}}

\newcommand{\setS}{\mathcal{S}}

\newcommand{\RomanOne}{\mathrm{I}}
\newcommand{\RomanTwo}{\mathrm{I\hspace{-.1em}I}}

\newcommand{\Dinv}{D^{-}}

\newtheorem{problem}{Problem}[section]

\newtheorem{algorithm}{Algorithm}[section]
\usepackage[hang,small]{caption}
\usepackage[subrefformat=parens]{subcaption}
\usepackage{algorithm}
\usepackage{algorithmicx}
\usepackage{algpseudocode}
\algdef{SE}[DOWHILE]{Do}{doWhile}{\algorithmicdo}[1]{\algorithmicwhile\ #1}

\newcommand{\sa}{S_{\mathrm{a}}}
\newcommand{\se}{S_{\mathrm{e}}}
\newcommand{\cd}{C_{\RomanTwo}}
\newcommand{\cb}{C_{\RomanOne}}
\newcommand{\indexc}{\mathfrak{I}_c}
\newcommand{\indexr}{\mathfrak{I}_r}
\newcommand{\diag}{\operatorname{diag}}
\usepackage{amsbsy}
\usepackage{color}
\usepackage{enumerate}

\usepackage{bm}

\usepackage{yhmath}
\usepackage{color}

\usepackage{arydshln}

\makeatletter

\def\env@sqcases{%
  \let\@ifnextchar\new@ifnextchar
  \left\lfloor\,\,\,
  \def\arraystretch{1.2}%
  \array{@{}l@{\quad}l@{}}%
}
\makeatother

\DeclareMathAlphabet\mathbfcal{OMS}{cmsy}{b}{n}

\DeclareFontFamily{U}{mathx}{\hyphenchar\font45}
\DeclareFontShape{U}{mathx}{m}{n}{
      <5> <6> <7> <8> <9> <10>
      <10.95> <12> <14.4> <17.28> <20.74> <24.88>
      mathx10
      }{}
\DeclareSymbolFont{mathx}{U}{mathx}{m}{n}
\DeclareMathSymbol{\bigtimes}{1}{mathx}{"91}

\begin{document}
\title{Imposing early and asymptotic constraints on LiGME with application\\ to bivariate nonconvex enhancement of fused lasso models}

\name{Wataru Yata and Isao Yamada\thanks{
  This work was supported partially by JSPS Grants-in-Aid (23KJ0945, 19H04134) and by JST SICORP (JPMJSC20C6).}}
\address{
  \textit{
  Department of Information and Communications Engineering}\\
  \textit{Tokyo Institute of Technology, 
Tokyo, Japan}\\
\{yata,isao\}@sp.ict.e.titech.ac.jp
}
\maketitle
\setlength\dashlinedash{.5pt}
\setlength\dashlinegap{1.3pt}

\abovedisplayskip= 2.5pt
\belowdisplayskip = 2.5pt
\begin{abstract}
For the constrained LiGME model, a nonconvexly regularized least squares estimation model, we present an iterative algorithm of guaranteed convergence to its globally optimal solution.
The proposed algorithm can deal with two different types of constraints simultaneously. 
The first type constraint, called the asymptotic one, requires the limit of estimation sequence to achieve the corresponding condition.
The second type constraint, called the early one, requires every vector in estimation sequence to achieve the condition. 
We also propose a bivariate nonconvex enhancement of fused lasso models with effective constraint
for sparse piecewise constant signal estimations.
(This is an improved version of [Yata and Yamada, ICASSP 2024].)
\end{abstract}

\section{Introduction}
In this paper, 
we revisit the following
signal estimation model (see Section \ref{sec:ea-type} for detail):
\begin{align}
  \label{eq:general-model}
  \minimize_{
    \substack{
    x\in C_0\\
    (j \in \indexc)\, \opC_j x\in C_j} 
  } J(x):=\frac{1}{2}
  \norm{y - A x}^2_{\spY} +\mu \sum_{i\in \indexr} \mu_i \Psi^{\lrangle{i}}_{B^{\lrangle{i}}}\circ \opL_i (x),
\end{align}
where $y$ is an observed signal
in a finite dimensional real Hilbert space $\spY$, $A$ is a linear operator, $\mu>0$, for each $i\in \indexr$, $\opL_i$ is a linear operator, $\mu_i > 0$ and $\Psi^{\lrangle{i}}_{B^{\lrangle{i}}}$ 
is a generalized Moreau enhancement of a prox-friendly convex function $\Psi^{\lrangle{i}}$ with a tuning matrix $B^{\lrangle{i}}$ (see Section \ref{sec:notation} and \ref{sec:pre-cLiGME}).
In \eqref{eq:general-model}, $C_0$ and $(j\in\indexc) \,C_j$ are simple closed convex sets meaning that the metric projections $P_{C_0}$ and $(j\in\indexc) \,P_{C_j}$ are available as computable operators, and $(j\in\indexc) \,\opC_j$ are linear operators.
The cost function $J$ in \eqref{eq:general-model} was introduced originally in the Linearly involved Generalized Moreau Enhanced (LiGME) model \cite{LiGME}, where each regularizer $(i\in \indexr)\, \Psi^{\lrangle{i}}_{B^{\lrangle{i}}}\circ \opL_i(x)$ is a linearly involved generalization
of the generalized minimax concave (GMC) penalty \cite{selesnick2017}
(and the MC penalty \cite{zhang2010}) and serves as a parametric bridge of the gap
between the direct discrete measures (e.g., $\norm{\cdot}_0$ and $\operatorname{rank}(\cdot)$) and their convex envelope functions (e.g., $\norm{\cdot}_1$ and $\norm{\cdot}_{\mathrm{nuc}}$) (see \cite[Example 2]{LiGME}).
With some proper tunings of $B^{\lrangle{i}}$ \cite{LiGME,liu2022,Chen2023}, we can achieve simultaneously the nonconvexity of $(i\in \indexr)\, \Psi^{\lrangle{i}}_{B^{\lrangle{i}}}$ and the overall convexity of $J$.
Moreover, under the overall convexity of $J$,
iterative algorithms were established for the LiGME model \cite{LiGME} and the constrained LiGME (cLiGME) model \cite{yata2021,cLiGME} in a frame of the fixed point strategies \cite{combettes2021fixedpointstra}
(for recent advancements of convexity
preserving nonconvex models, 
see, e.g., \cite{selesnick2017,LiGME,shabili2021,zhang2023}).

The cLiGME model \cite{yata2021, cLiGME} 
introduced \textit{the split feasibility type convex constraints}\footnote{
  Thanks to the great expressive ability of the split feasibility type constraint set, $\sa$ has been utilized as key information in many signal processing applications \cite{herman1980,sezan1992,byrne2004,censor2005}.
} \cite{byrne2004,censor2005}
$(j \in \indexc)\, \opC_j x\in C_j$ under the consistent condition $ \sa:= \bigcap_{j \in \indexc\,}\{ x\,| \opC_j x\in C_j\} \neq \emptyset$.
For the cLiGME model, \cite{yata2021,cLiGME} proposed an iterative algorithm that can produce a sequence $(x_k)_{k\in \setN}$ which is guaranteed to converge to a minimizer of $J$ over $\sa$.
In this sense, we call $\sa$ \textit{the asymptotic constraint set}\footnote{
  Even for
  \textit{
  the split feasibility problem} just to find one anonymous point in $\sa$, it is known that its direct computation is hard and iterative algorithms
  have been studied 
  for producing a sequence to converge to
  a point in $\sa$ \cite{byrne2004,censor2005}.
}.

On the other hand, 
in a medical image application \cite{kitahara2021} of the LiGME model \cite{LiGME}, a simple nonnegativity constraint $x \in C_0$ was introduced
to enforce $(\forall k \in \setN)\, x_k \in C_0$.
To distinguish from the asymptotic constraint set $\sa$,
we call such a type of constraint set applicable to enforce $(\forall k \in \setN)\, x_k \in C_0=:\se$ 
\textit{the early constraint set}.

In this paper, for further flexible applicability of the cLiGME model, we impose both the asymptotic constraint $\sa$ and
the early constraint $\se$
on the model \eqref{eq:general-model}
under $\sa\cap\se \neq \emptyset$
and present a modified iterative algorithm to solve problem \eqref{eq:general-model} under early and asymptotic constraints simultaneously.
The proposed algorithm (EA-type constrained algorithm) in Section \ref{sec:ea-type} can be seen as a unified extension of
algorithms in \cite{cLiGME,kitahara2021}.
The proposed algorithm produces a sequence $(x_k)_{k\in\setN}\subset C_0=\se$ of guaranteed convergence to a minimizer of
$J$ over $\sa\cap\se$
implying thus $x_{k_0}\in C_0=\se$ is ensured even in a case where we have to use $x_{k_0}$ for some finite $k_0$ before convergence of the algorithm due to the time limitations.

To verify the effectiveness of the cLiGME model with the proposed EA-type constrained algorithm, in Section \ref{sec:latent-fused-cLiGME}, 
we present an application to nonconvex and robustness enhancements of the fused lasso models \cite{tibshirani2005,feng2020latent} used 
for sparse piecewise constant signal estimations\footnote{
  Indeed, for piecewise continuous signal estimations, including piecewise constant signal estimations, certain hidden sparsities in the target information have been exploited successfully. See, e.g., \cite{angelosante2010,lin2017,kuroda2018}.
} \cite{tibishirani2007,mishra2023} possibly under nonzero baseline assumptions.

\section{Preliminaries}
\label{sec:Preliminaries}
\subsection{Notation}
\label{sec:notation}
Let $\spH$ and $\spK$ be finite dimensional real Hilbert spaces.
A Hilbert space $\spH$ is equipped with an inner product $\ip{\cdot}{\cdot}_\spH$
and its induced norm $\norm{\cdot}_\spH$. 
$\setLO{\spH}{\spK}$ denotes the set of all linear operators from 
$\spH$ to $\spK$.
For $L\in \setLO{\spH}{\spK}$, 
$\norm{L}_{\mathrm{op}}$ denotes the operator norm of $L$
(i.e., $\norm{L}_{\mathrm{op}}
:=\sup_{x\in\spH, \norm{x}_{\spH}\leq 1} \norm{Lx}_\spK$) and
$L^*\in\setLO{\spK}{\spH}$ the adjoint operator of $L$ 
(i.e., $(\forall x\in \spH)(\forall y \in\spK)\, 
\ip{Lx}{y}_\spK = \ip{x}{L^*y}_\spH$).
The identity operator is denoted by
$\Id$ and the zero operator from $\spH$ to $\spH$ by $\zeroMatrix_\spH$.
We express the positive definiteness and the positive semidefiniteness of
a self-adjoint operator $L\in \setLO{\spH}{\spH}$ 
as $L\succ \zeroMatrix_\spH$ and $L\succeq  \zeroMatrix_\spH$, respectively.
The set of all proper lower semicontinuous convex functions\footnote{A function $f:\spH\to (-\infty,\infty]$ is 
(i) proper if $\dom f\neq \emptyset$, 
(ii) lower semicontinuous
if $ \{ x\in \spH | f(x)\leq \alpha\}$ 
is closed for $\forall \alpha\in \setR$,
(iii) convex if 
$f(\theta x + (1-\theta) y)\leq \theta f(x)+(1-\theta) f(y) $
for $\forall x,y \in \spH,\, 0<\theta< 1$.} defined on $\spH$ is denoted by $\Gamma_0 (\spH)$.
$f\in\Gamma_0(\spH)$ is said to be prox-friendly if 
$\prox_{\gamma f}:\spH\to\spH:x\mapsto \argmin_{v\in\spH}\left[
  \Psi(v) + \frac{1}{2\gamma}\norm{v-x}^2_\spH
\right]$
is available as a computable operator for $\forall\gamma>0$.
A closed convex set $C\subset \spH$ is simple if $P_C:\spH\to\spH:x\mapsto \argmin_{y\in C}\norm{x-y}_\spH$ is available as a computable operator.
Note that the proximity operator of the indicator function\footnote{
  For a closed convex set $C$, $\iota_C$ is defined
as $
  \iota_C(x) := \begin{cases}
  0 & (x\in C)\\
  +\infty & (x\notin C).
\end{cases}$
} $\iota_C$ reproduces $P_C$.
$\mathbf{0}_n\in\setR^n$ ($\zeroMatrix_{m,n}\in \setR^{m\times n}$ or $\zeroMatrix_n\in \setR^{n\times n}$) denotes the zero vector (zero matrix). 
$\mathbf{1}_n\in\setR^n$ denotes the
all one vector.
For $A \in \setR^{m\times n}$, 
$A^\transT$ denotes the transpose of $A$.
$I_n\in \setR^{n\times n}$ stands for the identity matrix.
$\diag(A_1,A_2):= \begin{bmatrix}
  A_1 & \zeroMatrix_{m_2, n_1}\\
  \zeroMatrix_{m_1, n_2} & A_2
\end{bmatrix}\in\setR^{(m_1+m_2)\times (n_1+n_2)}$ denotes the block diagonal matrix of $A_1\in \setR^{m_1\times n_1}$ and $A_2\in \setR^{m_2\times n_2}$.
\subsection{cLiGME model and its overall convexity condition}
\label{sec:pre-cLiGME}
\begin{problem}[cLiGME model {\cite[Problem 3.1]{cLiGME}}]
  \label{prob:cLiGME}
	\thickmuskip=0\thickmuskip
	\medmuskip=0\medmuskip
	\thinmuskip=0.0\thinmuskip
	\arraycolsep=0.5\arraycolsep
  \sloppy
  Let 
  $\indexr$ and $\indexc$ be finite index sets. Let $\mathcal{X}$,
  $\mathcal{Y}$,
  $(i\in \indexr)\, \spZ_i, \sptildeZ_i$ and
  $(j\in \indexc)\, \mathfrak{Z}_j$
  be finite dimensional real Hilbert spaces. Suppose that  
  (a) $A\in \setLO{\spX}{\spY}$, $y\in\spY$ and $\mu>0$;
  (b) $(\forall i\in \indexr)\, B^{\lrangle{i}}\in \setLO{\spZ_i}{\sptildeZ_i}$, $\opL_i\in \setLO{\spX}{\spZ_i}$ and $\mu_i >0$;
  (c)\footnote{The condition (c) can be relaxed, e.g., $\Psi^{\lrangle{i}}$ introduced in \eqref{eq:latent-fused-cLiGME} violates $\dom \Psi^{\lrangle{i}}= \spZ_i$ but can  be used in models \eqref{eq:original-cLiGME-model} and \eqref{eq:ea-cLiGME} and in Algorithm \ref{alg:enforced}.}  
  $(\forall i\in \indexr)\,\Psi^{\lrangle{i}}\in \Gamma_0(\spZ_i)$ is
  (i) coercive, (ii) $\dom\Psi^{\lrangle{i}} = \spZ_i$,
  (iii) even symmetry (i.e., $\Psi^{\lrangle{i}}\circ(-\Id) =\Psi^{\lrangle{i}}$),
  (iv) prox-friendly;
  (d) $(\forall j\in \indexc)\, C_j(\subset\spfrakZ_j)$ is a simple closed convex set and $\opC_j\in \setLO{\spX}{\spfrakZ_j}$.
  Assume $\bigcap_{j \in \indexc\,}\{ x\in \spX\,| \opC_j x\in C_j\} \neq \emptyset$.
  Then 
  \begin{enumerate}[(a)]
    \item For each $\Psi^{\lrangle{i}}$, 
    its generalized Moreau enhancement is defined by 
    $
      \Psi^{\lrangle{i}}_{B^{\lrangle{i}}}(\cdot):= \Psi^{\lrangle{i}}(\cdot) - \min_{v\in\spZ_i}\left[
        \Psi^{\lrangle{i}}(v) +\frac{1}{2}\norm{B^{\lrangle{i}}(\cdot-v)}^2_{\sptildeZ_i}
      \right]
    $
    with a tuning operator $B^{\lrangle{i}}\in \setLO{\spZ_i}{\sptildeZ_i}$ called GME matrix.
    \item The constrained LiGME (cLiGME) model is formulated as
    \begin{align}
      \label{eq:original-cLiGME-model}
      \minimize_{(j \in \indexc)\, \opC_j x\in C_j} 
      J(x) \,\,\,\,\quad (\mbox{see \eqref{eq:general-model}}).
    \end{align}
  \end{enumerate}
\end{problem}

Although the model \eqref{eq:original-cLiGME-model} is seemingly different from the original cLiGME model \cite[Problem 3.1]{cLiGME}, the model \eqref{eq:original-cLiGME-model} can be translated into \cite[Problem 3.1]{cLiGME} via a product space reformulation \cite[Remark 3.1, Remark 3.2]{cLiGME}.
Note that, even if $\Psi^{\lrangle{1}}$ and $\Psi^{\lrangle{2}}$ are prox-friendly, $\Psi^{\lrangle{1}}\circ \opL_1 + \Psi^{\lrangle{2}}$ is not necessarily prox-friendly.
Therefore, we only assume that $(i\in \indexr)\, \Psi^{\lrangle{i}}$ are all prox-friendly.

\begin{example}[Selected reproductions from \eqref{eq:original-cLiGME-model}]
  \label{example:reproduction}
  \quad
  \begin{enumerate}[(a)]
    \item (Fused lasso \cite{tibshirani2005}). 
    For estimation of a sparse piecewise constant signal $\xstar\in\setR^N$ from its noisy observation $y (\approx A\xstar) \in\setR^M$
    with $A\in\setR^{M\times N}$, the fused lasso model estimates $\xstar$ by
    \begin{align}
      \label{eq:fused-lasso}
      \widetilde{x}\in
    \argmin_{x\in\setR^N} \frac{1}{2}\norm{A x -y}^2 + \mu_1 \norm{x}_1 +\mu_2 \norm{D x}_1,
    \end{align}
    where $\mu_1,\mu_2\geq 0$ are chosen regularization weights (the lasso model \cite{tibshirani1996} is reproduced simply by $\mu_2=0$) and 
    \begin{align*}
      D := \begin{bmatrix}
        -1 & 1      &        & \\
        & \ddots & \ddots  & \\
        &        & -1      & 1
      \end{bmatrix}\in \setR^{(N-1)\times N}.
    \end{align*}
    The model \eqref{eq:fused-lasso} is an instance of the model \eqref{eq:original-cLiGME-model}. 

    \item (GMC model \cite{selesnick2017}). By $|\indexr|:=1$, $\indexc:=\emptyset$, $\spX:= \setR^N$,$\spY:= \setR^M$, $\spZ:=\setR^N$, 
    $\sptildeZ:=\setR^M$,
    $\Psi := \norm{\cdot}_1$ and $\opL:= I_N$,
    the model \eqref{eq:original-cLiGME-model} reproduces the generalized minimax concave model.
  \end{enumerate}
\end{example}

The following fact shows the overall convexity condition for the cLiGME model \eqref{eq:original-cLiGME-model}.
\begin{proposition}[Overall convexity condition for the model \eqref{eq:original-cLiGME-model}
  {\cite[Proposition 1]{LiGME}, \cite[Proposition 3.2]{cLiGME}}]
  \label{prop: overall convexity}
  \sloppy
    Consider the settings in Problem \ref{prob:cLiGME}.
    For three
    conditions
    $(C_1)\ A^*A - \mu\sum_{i\in \indexr} \mu_i\opL_i^*{B^{\lrangle{i}}}^*{B^{\lrangle{i}}}\opL_i
    \succeq \zeroMatrix_\mathcal{X}$,
    $(C_2)\ J\in\Gamma_0(\mathcal{X})
    \mbox{ for any }y\in \mathcal{Y}$,
    $(C_3)\ J^{(0)}:=
    \frac{1}{2}\norm{A\cdot}^2_\mathcal{Y}+
    \mu\sum_{i\in \indexr} \mu_i \Psi^{\lrangle{i}}_{B^{\lrangle{i}}}\circ \opL_i\in\Gamma_0(\mathcal{X})$, 
    the relation $(C_1)\implies (C_2)\iff (C_3)$ holds.
    We call the condition $(C_1)$ the overall convexity condition.
\end{proposition}
For designs of $(i\in \indexr)\, B^{\lrangle{i}}$
achieving the overall convexity condition $(C_1)$, see \cite{LiGME,liu2022,Chen2023}.

\subsection{Latent fused lasso}
\label{sec:fused-lasso}
In the fused lasso \eqref{eq:fused-lasso}, $\norm{x}_1$ and $\norm{Dx}_1$ are used in the regularization terms to promote the sparsity of $x$ and the piecewise constancy of $x\in\setR^N$.
However, the gap between $\norm{\cdot}_0$ and $\norm{\cdot}_1$ causes underestimation of the estimated signal values by the fused lasso \eqref{eq:fused-lasso}.
To reduce such doubly underestimation effect, the latent fused lasso \cite{feng2020latent} was designed to employ $\norm{\cdot}_1$ only to promote the piecewise constancy of $x$ while the sparsity of $x$ is promoted by a newly introduced constraint on synthesis components $(1\leq j \leq N-L)\, d_j\in \setR^L$
of the first order difference signal $Dx$, where $L$ is the length of synthesis components.
Indeed, the latent fused lasso \cite{feng2020latent} translates the short duration of the piecewise constant parts observed in the target signal into constraints $(1\leq j \leq N-L)\, \mathbf{1}^\transT_L d_j = 0$ (see also Problem \ref{prob:unification-and-robustification}).

Clearly, another possible advancement of the latent fused lasso \cite{feng2020latent} would be thoughtful application of nonconvex enhancement of $\norm{\cdot}_1$ in \cite[(13)]{feng2020latent} (see also \cite[Section 6]{feng2020latent}). 
In Section \ref{sec:latent-fused-cLiGME}, we propose such a nonconvex enhancement of the latent fused lasso model.

\section{EA-type constrained algorithm for cLiGME}
\label{sec:ea-type}
First, we describe a detail of the model \eqref{eq:general-model}.
\begin{problem}[EA-type constrained LiGME]
Let
the settings in Problem \ref{prob:cLiGME} be in force. Let $C_0\subset \spX$ be a simple closed convex set.
Assume $C_0\cap\left(\bigcap_{j \in \indexc\,}\{ x\in\spX\,| \opC_j x\in C_j\}\right)\neq \emptyset$.
Then consider
\begin{align}
  \label{eq:ea-cLiGME}
  \minimize_{\substack{x\in C_0\\
  (j \in \indexc)\, \opC_j x\in C_j}}
  J(x) \,\,\,\,\quad (\mbox{see \eqref{eq:general-model}}).
\end{align}
\end{problem}
Note that, although the model \eqref{eq:ea-cLiGME} can be reformulated to the cLiGME model \eqref{eq:original-cLiGME-model} by a product space reformulation \cite[Remark 3.1]{cLiGME}, we distinguish $C_0$ as the early constraint set $\se$ from
the asymptotic constraint set $\sa$.
To minimize the model \eqref{eq:ea-cLiGME}
with guaranteed satisfaction of the early constraint condition
$(k\in \setN)\, x_k\in C_0$, 
we propose a EA-type algorithm for the model \eqref{eq:ea-cLiGME}
(Algorithm \ref{alg:enforced}).
\begin{algorithm}[h]
	\thickmuskip=0\thickmuskip
	\medmuskip=0\medmuskip
	\thinmuskip=0.0\thinmuskip
	\arraycolsep=0.5\arraycolsep
  \caption{EA-type constrained algorithm for the model \eqref{eq:ea-cLiGME}}
  \label{alg:enforced}
  \begin{algorithmic}[1]
    \State Choose $h_0:=\left(x_0,(v_0^{\lrangle{i}})_{i\in \indexr}, (w_0^{\lrangle{i}})_{i\in \indexr}, (z_0^{\lrangle{j}})_{j\in \indexc}\right)\in \spH:=\spX \times \left(\bigtimes_{i\in \indexr} \spZ_i\right) \times \left(\bigtimes_{i\in \indexr} \spZ_i\right) \times \left(\bigtimes_{j\in \indexc} \spfrakZ_j\right)$ satisfying $x_0\in C_0$.
    \State Choose $(\sigma,\tau,\kappa)
    \in (0,\infty)\times (0,\infty)\times(1,\infty)$ 
    satisfying\footnotemark
    \vspace*{-0.5em}
    \begin{align*}
      \begin{cases}
        \sigma\Id - \frac{\kappa}{2}A^*A- \mu
        \left(\sum_{i\in \indexr}  \opL_i^*\opL_i 
        + \sum_{j\in \indexc}\opC_j^*\opC_j\right)\succ \zeroMatrix_\spX \\
    \tau\geq \mu\left(\frac{\kappa}{2}+ \frac{2}{\kappa}\right)\max\left\{\mu_i\norm{B^{\lrangle{i}}}^2_{\mathrm{op}} \middle| i\in \indexr\right\}.
      \end{cases}
    \end{align*}
    \State Define  $T_{\mathrm{EA}}:\spH\to\spH:\left(x,(v^{\lrangle{i}})_{i\in \indexr}, (w^{\lrangle{i}})_{i\in \indexr}, (z^{\lrangle{j}})_{j\in \indexc}\right)\mapsto \left(\xi,(\zeta^{\lrangle{i}})_{i\in \indexr}, (\eta^{\lrangle{i}})_{i\in \indexr}, (\varsigma^{\lrangle{j}})_{j\in \indexc}\right)$ by 
    \Statex $\xi :=
      P_{C_0}\left[\left(\Id - 
      \frac{1}{\sigma}\left(A^*A 
      - \mu\sum_{i\in \indexr}\mu_i
      \opL_i^*{B^{\lrangle{i}}}^*B^{\lrangle{i}}\opL_i\right)
      \right)x+ \frac{1}{\sigma}A^*y
      \right.$
    \Statex \quad\quad\quad
      $
      -\frac{\mu}{\sigma}\sum_{i\in \indexr}
      \left(\mu_i\opL_i^*{B^{\lrangle{i}}}^*B^{\lrangle{i}}v^{\lrangle{i}}
       + \opL_i^*w^{\lrangle{i}}\right)\left. - 
       \frac{\mu}{\sigma}\sum_{j\in \indexc} \opC_j^* z^{\lrangle{j}}\right],$
      \Statex
      $\zeta^{\lrangle{i}}
      :=\prox_{\frac{\mu\mu_i}{\tau}\Psi^{\lrangle{i}}}\left[
        \frac{\mu\mu_i}{\tau}{B^{\lrangle{i}}}^*B^{\lrangle{i}} \opL_i (2\xi -x)
        +\left(\Id - \frac{\mu\mu_i}{\tau}{B^{\lrangle{i}}}^*B^{\lrangle{i}}\right)v^{\lrangle{i}}
        \right]$,
      \Statex
      $\eta^{\lrangle{i}}:=
        (\Id - \prox_{\mu_i\Psi^{\lrangle{i}}})\left(2\opL_i \xi-
        \opL_i x+w^{\lrangle{i}}\right)$,
        \Statex
        $\varsigma^{\lrangle{j}} :=
        (\Id -P_{C_j})\left(2\opC_j \xi-
        \opC_j x+ z^{\lrangle{j}}\right)$.

    \For{$k=0,1,2,\cdots$ }
    \State $h_{k+1} = T_{\mathrm{EA}}(h_k)$.
    \EndFor
  \end{algorithmic}
\end{algorithm}
\footnotetext{
	\thickmuskip=0\thickmuskip
	\medmuskip=0\medmuskip
	\thinmuskip=0.0\thinmuskip
	\arraycolsep=0.5\arraycolsep
  \sloppy For example, choose $\kappa>1$ and compute $(\tau,\kappa)$ by 
  \vspace*{-0.5em}
\begin{align*}
\begin{cases}
  \sigma :=\norm{\frac{\kappa}{2}A^*A +\mu \sum_{i\in \indexr} \opL_i^*\opL_i +\mu \sum_{j\in \indexc} \opC_j^*\opC_j}_{\mathrm{op}} + (\kappa-1)\\
  \tau:= \mu\left(\frac{\kappa}{2} + \frac{2}{\kappa}\right)\max\left\{\mu_i\norm{B^{\lrangle{i}}}^2_{\mathrm{op}} \middle| i\in \indexr\right\} + (\kappa -1).
\end{cases}
\end{align*}
}

Algorithm \ref{alg:enforced} is realized based on integration of
\cite[Algorithm 1]{cLiGME} and \cite[(16)]{kitahara2021}.
If $P_{C_0}= \Id$, Algorithm \ref{alg:enforced} reproduces \cite[Algorithm 1]{cLiGME}.
On the other hand, if $\indexc = \emptyset$, Algorithm \ref{alg:enforced} reproduces \cite[(16)]{kitahara2021}.
Under mild assumptions (see \cite[Assumption 3.1]{cLiGME} for detail), we provide (i) a guarantee of the early constraint condition and (ii) a convergence analysis.
\begin{theorem}
  \label{thm:convergence}
  \sloppy
  Consider the model \eqref{eq:ea-cLiGME} under the overall convexity condition. 
  Generate the sequence 
  $\left(h_k\right)_{k\in\setN}$
  by Algorithm \ref{alg:enforced} and set the sequence $(x_k)_{k\in \setN}$ by $(\forall k\in\setN)\, x_k := \Xi(h_k)$, where $\Xi:\spH\to\spX:\left(x,(v^{\lrangle{i}})_{i\in \indexr}, (w^{\lrangle{i}})_{i\in \indexr}, (z^{\lrangle{j}})_{j\in \indexc}\right)\mapsto x $.
  Then the sequence $(x_k)_{k\in \setN}$ 
  (i) satisfies the early constraint condition $(\forall k\in \setN)\, x_k\in C_0$ and (ii) converges to a global minimizer of \eqref{eq:ea-cLiGME}, i.e.,
  \begin{align*}
    \lim_{k\to\infty} x_k\in \setS:= \argmin_{
    \substack{
    x\in C_0\\
    (j \in \indexc)\, \opC_j x\in C_j} 
  } J(x).
  \end{align*}
\end{theorem}
\noindent
\textbf{(Proof sketch)}
\sloppy
(i) Clear by the design of $T_{\mathrm{EA}}$ in Algorithm \ref{alg:enforced}.
(ii) The solution set of \eqref{eq:ea-cLiGME}  
coincides with $\Xi(\fix (T_{\mathrm{EA}}))$,
where $h\in\fix (T_{\mathrm{EA}})\iff h =T_{\mathrm{EA}}(h)$.
Moreover, $T_{\mathrm{EA}}$ is averaged nonexpansive with a certain norm and thus a fixed point of $T_{\mathrm{EA}}$ can be obtained by Krasnosel'ski\u{\i}-Mann iteration \cite[Theorem 5.14]{CAaMOTiH}.

\section{Application to Nonconvex and robustness enhancements of fused lasso models}
\label{sec:latent-fused-cLiGME}
\subsection{A robustness enhancement of fused lasso models}
\begin{figure}[h]
  \centering
    \includegraphics[bb = 0 0 1160 679, scale=0.1]{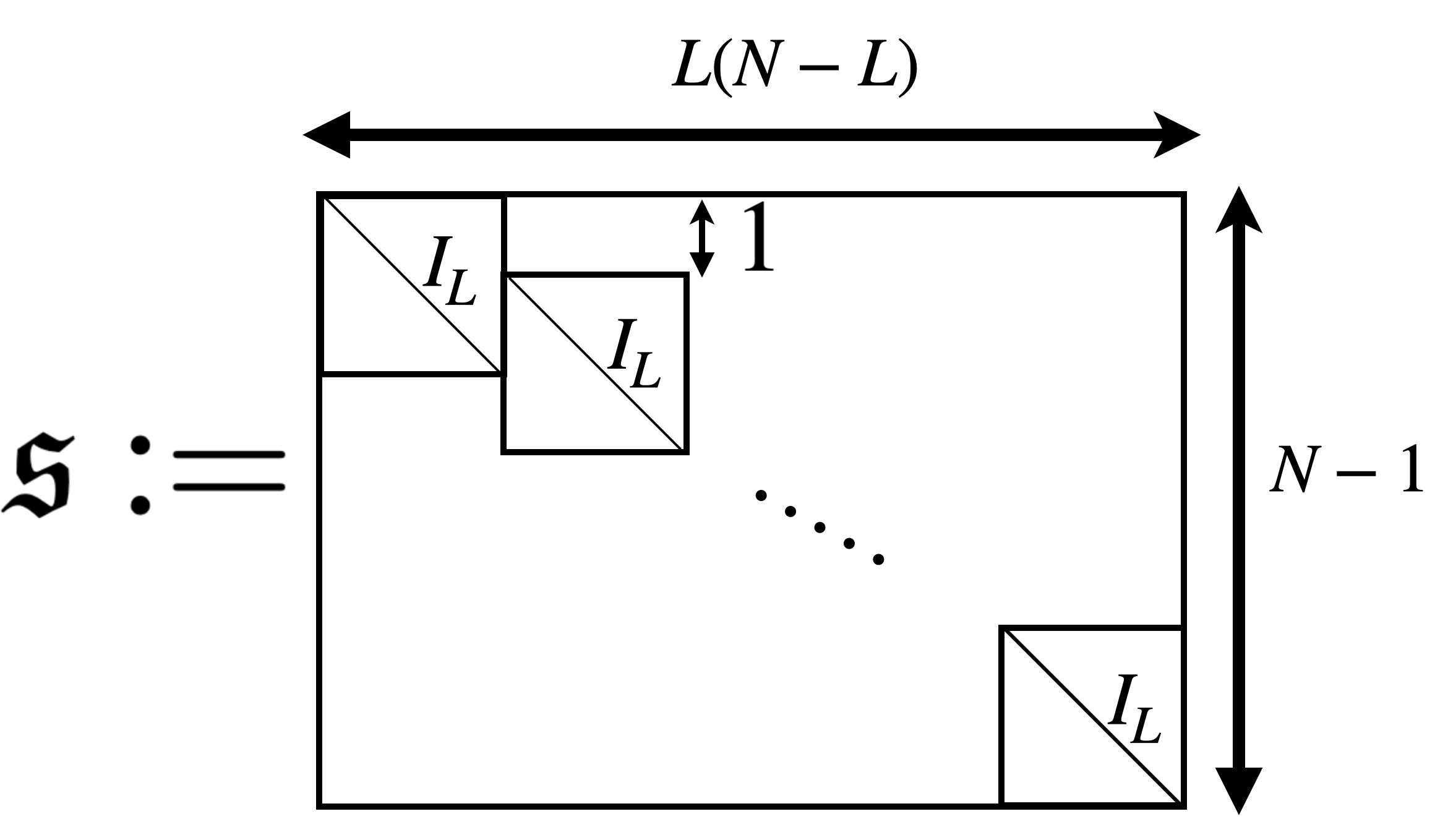}
  \caption{Illustration of linear operator $\mathfrak{s}$ by matrix form}
  \label{fig:def-S}
\end{figure}

\thickmuskip=0\thickmuskip
\medmuskip=0\medmuskip
\thinmuskip=0.0\thinmuskip
\arraycolsep=0.5\arraycolsep
 We reuse $A\in\setR^{M\times N}$, $y\in\setR^M$, $D\in \setR^{(N-1)\times N}$ according to Example \ref{example:reproduction}(a).
We also use, with  the length $L\in \{1,\ldots, N-1\}$ of synthesis components of $Dx$ (see Sec.\ref{sec:fused-lasso}), $\mathfrak{s}\in \setR^{(N-1)\times L(N-L)}$ (see Fig.\ref{fig:def-S}) and
\begin{align*}
  \Dinv = \begin{bmatrix}
    0 & \cdots      &        & 0\\
    1 & 0 &   & \\
    1 & 1 &    \ddots   &\vdots \\
    \vdots  &    &   \ddots    & 0\\ 
    1 & 1& \cdots& 1  
  \end{bmatrix}\in \setR^{N\times (N-1)},
  H= \begin{bmatrix}
    \mathbf{1}_{L}^\transT & & \\
     & \ddots & \\
     & & \mathbf{1}_{L}^\transT
  \end{bmatrix}\in \setR^{(N-L)\times L(N-L)}.
\end{align*}
\begin{problem}[Unification and robustification of fused lasso models]
  \label{prob:unification-and-robustification}
  Let $\mathcal{I}$ be a closed convex subset of $\setR$ and $\cb^{(\mathcal{I})}:= \{\alpha \mathbf{1}_N \mid \alpha\in\mathcal{I}\}$ a simple closed convex subset of $\cb^{(\setR)}:=\{ \alpha \mathbf{1}_N \mid \alpha \in \setR\}
  = \{ b\in\setR^N\mid Db= \mathbf{0}_{N-1}
  \} $, which is introduced to exploit a priori knowledge on the baseline vector in the sparse piecewise constant signal estimation.
  Let $(\mathbf{0}_{N-L}\in)\cd \subset \setR^{N-L}$ be a simple closed convex set, which is introduced possibly to enhance robustness against model mismatch in the latent fused lasso model \cite[(13)]{feng2020latent}. We propose the following equivalent models\footnote{Extensions of $\cb^{(\mathcal{I})}$ and $\cd$ will be discussed elsewhere.} (to be enhanced nonconvexly in Problem \ref{prob:ne-fused-lasso}).
	\thickmuskip=0\thickmuskip
	\medmuskip=0\medmuskip
	\thinmuskip=0.0\thinmuskip
	\arraycolsep=0.5\arraycolsep
  
  \begin{enumerate}
    \item[(a)] (A unified robustification of fused lasso models \eqref{eq:fused-lasso} and \cite[(13)]{feng2020latent}). Estimate $\xstar$ by $\widetilde{x}:= \widetilde{b} + \Dinv\mathfrak{s}(\widetilde{d})$ with
    \begin{align}
      \label{eq:unified-robustification}
      (\widetilde{b}, \widetilde{d})\in \argmin_{\substack{b\in\cb^{(\mathcal{I})}\\
        Hd\in\cd}} \frac{1}{2}\norm{y - A(b + \Dinv\mathfrak{s} (d))}^2 + \mu_1\norm{\Dinv\mathfrak{s}(d)}_1+ \mu_2\norm{d}_1, \nonumber \\[-25pt]
        \quad
    \end{align}
    where
  (i) $\mu_1,\mu_2\geq 0$ are chosen regularization weights,
  (ii) $b$ is an explicit baseline vector, and (iii) $d:= (d_1^\transT,\ldots, d_{N-L}^\transT)^\transT$ indicates synthesis components $d_j\in \setR^L\, (1\leq j \leq N-L)$ of $Dx$.

    \item[(b)] (Equivalent reformulations of \eqref{eq:unified-robustification}).
    By using 
    $\mathcal{A}:=A[\begin{array}{c:c} I_N & {\Dinv \mathfrak{s}} 
    \end{array}]\in\setR^{M\times (N+L(N-L))}$, \eqref{eq:unified-robustification} can be reformulated equivalently to:
    \begin{align}
      \label{eq:simplification}
      \minimize_{\substack{b\in\cb^{(\mathcal{I})}, Hd\in\cd\\ \mathbf{u}= (
        b^\transT, d^\transT
      )^\transT }}
      \frac{1}{2}\norm{y - \mathcal{A}\mathbf{u}
      }^2 + 
      \mu_1\norm{\Dinv\mathfrak{s}(d)}_1
      + \mu_2\norm{d}_1,
    \end{align}
    or 
    \begin{align}
      \label{eq:reformulation-2}
      \hspace{-1.5em}
      \minimize_{\substack{b\in\cb^{(\mathcal{I})},
      Hd\in\cd\\
      \mathbf{u}= (
        b^\transT, d^\transT
      )^\transT }}
      &\frac{1}{2}\norm{y - \mathcal{A}\mathbf{u}
      }^2 + 
      \mu_1(\iota_{\{\mathbf{0}_{N-1}\}} \oplus \norm{\cdot}^{(a)}_1) \circ \opL_1 (\mathbf{u}) \nonumber \\[-15pt]
      &+ \mu_2( \iota_{\{\mathbf{0}_{N-1}\}}\oplus \norm{\cdot}^{(b)}_1)\circ \opL_2 (\mathbf{u}),
    \end{align}
    where\footnote{
      $\norm{\cdot}^{(a)}_1$ 
      and $\norm{\cdot}^{(b)}_1 $ stand for $l_1$ norm functions defined respectively on $\setR^{N}$ and on $\setR^{L(N-L)}$.
    } $\opL_1 :=\diag(D, \Dinv\mathfrak{s})$, $\opL_2:= \diag(D, I_{L(N-L)})$ and, e.g., $(\iota_{\{\mathbf{0}_{N-1}\}} \oplus \norm{\cdot}^{(b)}_1) (r_1,r_2):=  \iota_{\{\mathbf{0}_{N-1}\}}(r_1) + \norm{r_2}^{(b)}_1 \,  (r_1\in\setR^{N-1},r_2\in\setR^{L(N-L)})$.
  \end{enumerate}
\end{problem}

\begin{remark}[Problem \ref{prob:unification-and-robustification}
  reproduces existing models]
  \label{remark:generalization-of-fls}
  \quad
  \begin{enumerate}[(a)]
  \item (Fused lasso \eqref{eq:fused-lasso}). By setting $L= N-1$, the third term in
  the cost function of \eqref{eq:unified-robustification} reproduces the TV regularizer: $\norm{\cdot}_1\circ D (x)= \norm{d}_1$ 
  (see \cite[Remark]{feng2020latent}). By $  \Dinv\mathfrak{s}(d) = x- b$,
  the model \eqref{eq:unified-robustification} can be seen as a constraint enhancement of the fused lasso \eqref{eq:fused-lasso}.
  \item (Latent fused lasso \cite{feng2020latent}). By setting 
  $M:=N$, $A := I_N$, $\mu_1:=0$, $\mu_2>0$,
  $\cd:= \{\mathbf{0}_{N-L}\}\subset \setR^{N-L}$
  and $\cb^{(\mathcal{I})} := \cb^{(\setR)}$, 
  \eqref{eq:unified-robustification}
  reproduces an equivalent model of the latent fused lasso \cite[(13)]{feng2020latent}.
  \end{enumerate}
\end{remark}

\subsection{Bivariate nonconvex enhancement of fused lass}
\begin{problem}[Bivariate nonconvex enhancement of \eqref{eq:reformulation-2}]
  \label{prob:ne-fused-lasso}
\thickmuskip=0\thickmuskip
\medmuskip=0\medmuskip
\thinmuskip=0.0\thinmuskip
\arraycolsep=0.5\arraycolsep
Estimate the target signal $\xstar$ by $\widetilde{x}:= \widetilde{b} + \Dinv\mathfrak{s}(\widetilde{d})$ with
\begin{align}
  \label{eq:latent-fused-cLiGME}
  (\widetilde{b}, \widetilde{d})\in\argmin_{\substack{
    b\in\cb^{(\mathcal{I})},
    Hd \in \cd\\
    \mathbf{u} = (
        b^\transT,
        d^\transT
      )^\transT 
  }} 
  &\frac{1}{2}\norm{y - \mathcal{A}\mathbf{u}
  }^2 + 
  \mu_1 (\iota_{\{\mathbf{0}_{N-1}\}}\oplus\norm{\cdot}^{(a)}_1)_{B^{\lrangle{1}}}\circ\opL_1(\mathbf{u}) \nonumber\\[-15pt]
  &\vspace{10em}+\mu_{2} (\iota_{\{\mathbf{0}_{N-1}\}} \oplus \norm{\cdot}^{(b)}_1)_{B^{\lrangle{2}}}\circ\opL_2 (\mathbf{u}),
\end{align}
where $B^{\lrangle{1}}\in \setR^{l_1\times (2N-1)}$ and $B^{\lrangle{2}}\in \setR^{l_2\times (N-1 +L(N-L))}$ are designed to satisfy the overall convexity condition (see, e.g., Proposition \ref{prop: overall convexity}).
\end{problem}

\begin{remark}[Why bivariate nonconvex enhancement of \eqref{eq:reformulation-2} ?]\quad
  \label{remark:why-nonconvex}
  \begin{enumerate}[(a)]
    \item
    Observing the following relation together with footnote 5:
  \thickmuskip=0\thickmuskip
  \medmuskip=0\medmuskip
  \thinmuskip=0.0\thinmuskip
  \arraycolsep=0.5\arraycolsep
    \begin{align*}
      &(\forall(b,d)\in \cb^{(\mathcal{I})}\times\setR^{L(N-L)})\ \ (\iota_{\{\mathbf{0}_{N-1}\}}\oplus \norm{\cdot}^{(b)}_1)_{B^{\lrangle{2}}}(Db, d)\\
      = &\norm{d}^{(b)}_1 
      -\min_{v_2\in \setR^{L(N-L)}}
      \norm{v_2}^{(b)}_1
      + \frac{1}{2}\norm{B^{\lrangle{2}}(\mathbf{0}_{N-1}^\transT, (d-v_2)^\transT)^\transT }^2\\
      =  &\norm{d}^{(b)}_1 
      -\min_{v_2\in \setR^{L(N-L)}}\left[
      \norm{v_2}^{(b)}_1
      + \frac{1}{2}\norm{B^{\lrangle{2}}_{ri}(d-v_2) }^2\right]
      = (\norm{\cdot}^{(b)})_{B^{\lrangle{2}}_{ri}}(d),
    \end{align*}
    we see that, under $b\in\cb^{(\mathcal{I})}$, $(\iota_{\{\mathbf{0}_{N-1}\}}\oplus \norm{\cdot}^{(b)}_1)_{B^{\lrangle{2}}}(Db,d)$ achieves a nonconvex enhancement of $\norm{d}^{(b)}_1$, with a GME matrix $B^{\lrangle{2}}_{ri}\in\setR^{l_2\times L(N-L)}$ ($B^{\lrangle{2}}_{ri}$ is the right block matrix of $B^{\lrangle{2}}$), without any influence on the baseline vector $b$. This observation can also be applied to interpretation of $(\iota_{\{\mathbf{0}_{N-1}\}}\oplus\norm{\cdot}^{(a)}_1)_{B^{\lrangle{1}}}(Db, \Dinv\mathfrak{s}(d))$ 
    as a nonconvex enhancement of $\norm{\Dinv\mathfrak{s}(d)}^{(a)}_1$.

    \item We found recently that (i) unfortunately, the design of $(B^{\lrangle{1}},B^{\lrangle{2}})$ reported in \cite{yata2024} does not really guarantee the overall convexity over the full space $\setR^{N}\times \setR^{L(N-L)}$
    essentially because the overall convexity condition over $\setR^{N}\times \setR^{L(N-L)}$ is too restrictive 
    in this scenario, (ii) we can achieve the bivariate overall convexity by designing $(B^{\lrangle{1}},B^{\lrangle{2}})$ with use of constraints. 
    \item Fortunately, we can design effective bivariate GME matrices $(B^{\lrangle{1}}, B^{\lrangle{2}})$ in \eqref{eq:latent-fused-cLiGME} by applying \cite[Corollary 1]{LiGME} together with\footnote{The GME matrix design in \cite[Theorem 1]{Chen2023} is made with LDU decomposition of $\opL_1$ and $\opL_2$ in \eqref{eq:latent-fused-cLiGME}.} \cite[Theorem 1]{Chen2023} in place of \cite[Proposition 2]{LiGME}.
\end{enumerate}
\end{remark}

For \eqref{eq:latent-fused-cLiGME},
we propose to use Algorithm \ref{alg:enforced} with the early constraint
$b\in \cb^{(\mathcal{I})}$ and the asymptotic constraint\footnote{
  In general, $\{d \in \setR^{L(N-L)} |\, H d\in \cd \}$
  is not simple even if $\cd$ is simple.
} $Hd\in \cd$.

\subsection{Numerical experiment}
\begin{figure}[H]
  \centering
  \begin{minipage}[b]{0.46\linewidth}
    \centering
    \includegraphics[bb = 0 0 1120 840, scale=0.085]
    {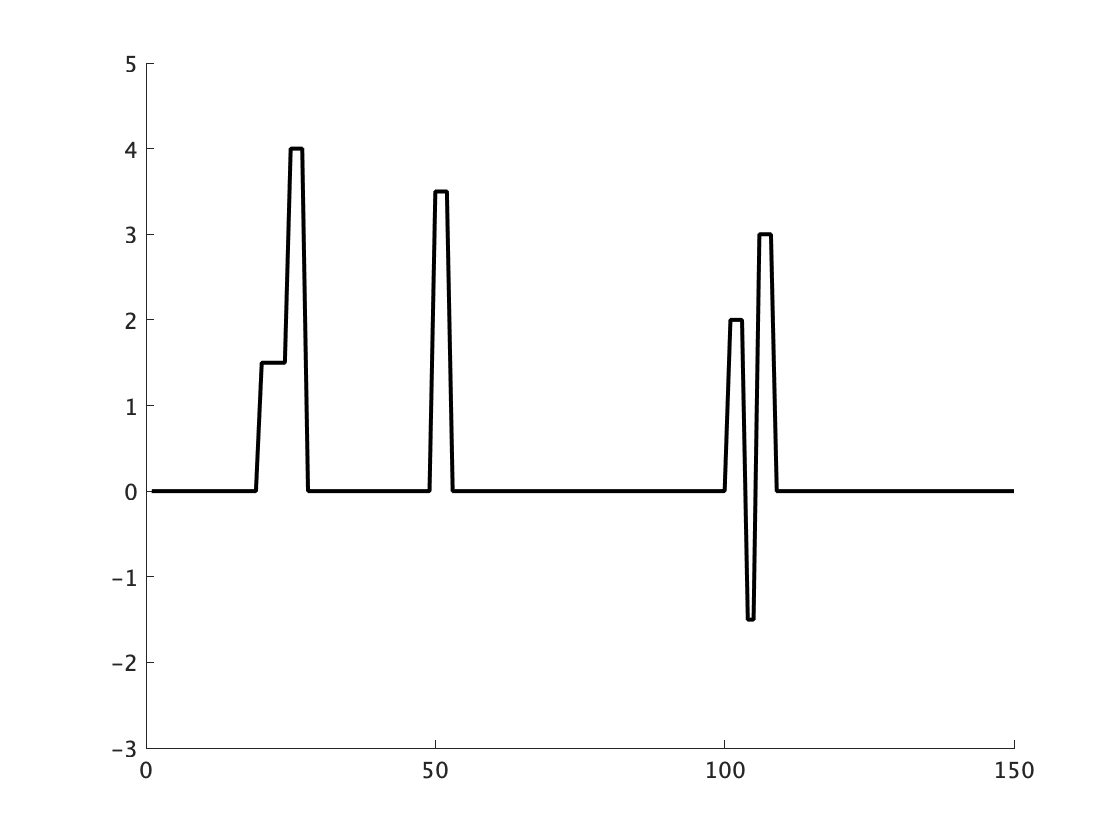}
    \subcaption{}
  \end{minipage}
  \begin{minipage}[b]{0.46\linewidth}
    \centering
    \includegraphics[bb = 0 0 1120 840, scale=0.085]
    {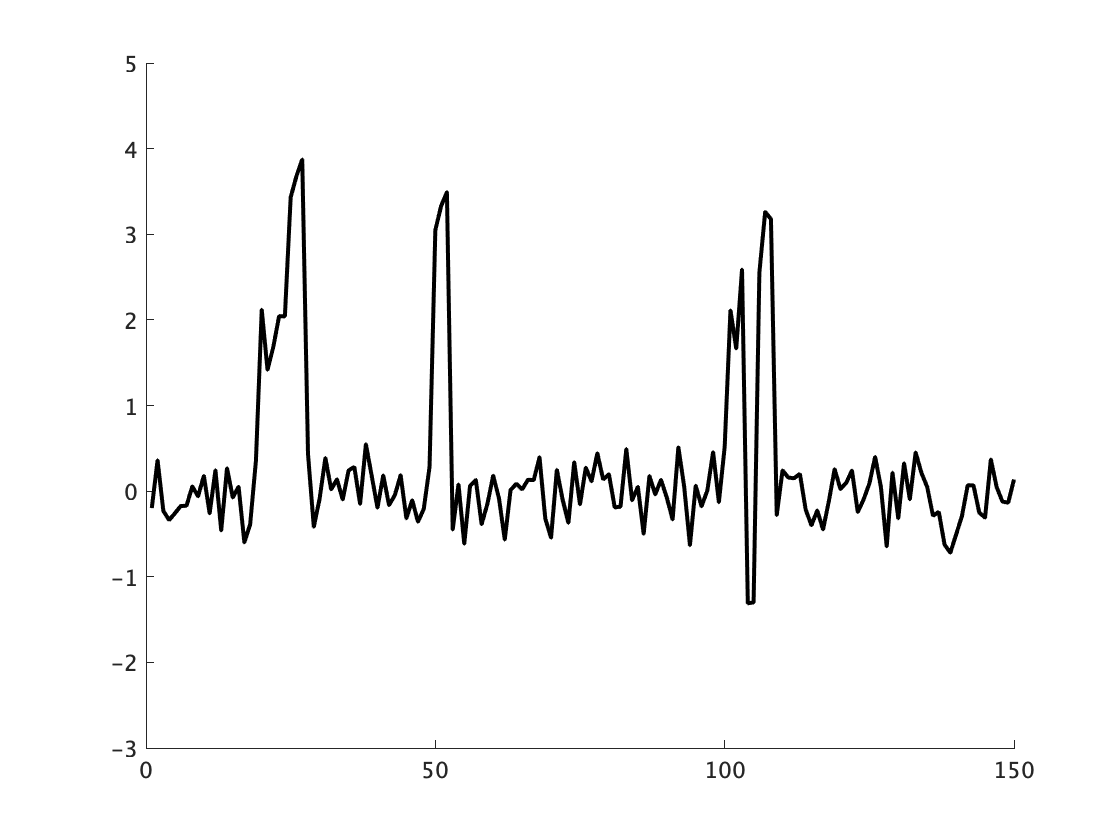}
    \subcaption{}
  \end{minipage}\\
  \begin{minipage}[b]{0.46\linewidth}
    \centering
    \includegraphics[bb = 0 0 1120 840,scale=0.085]
    {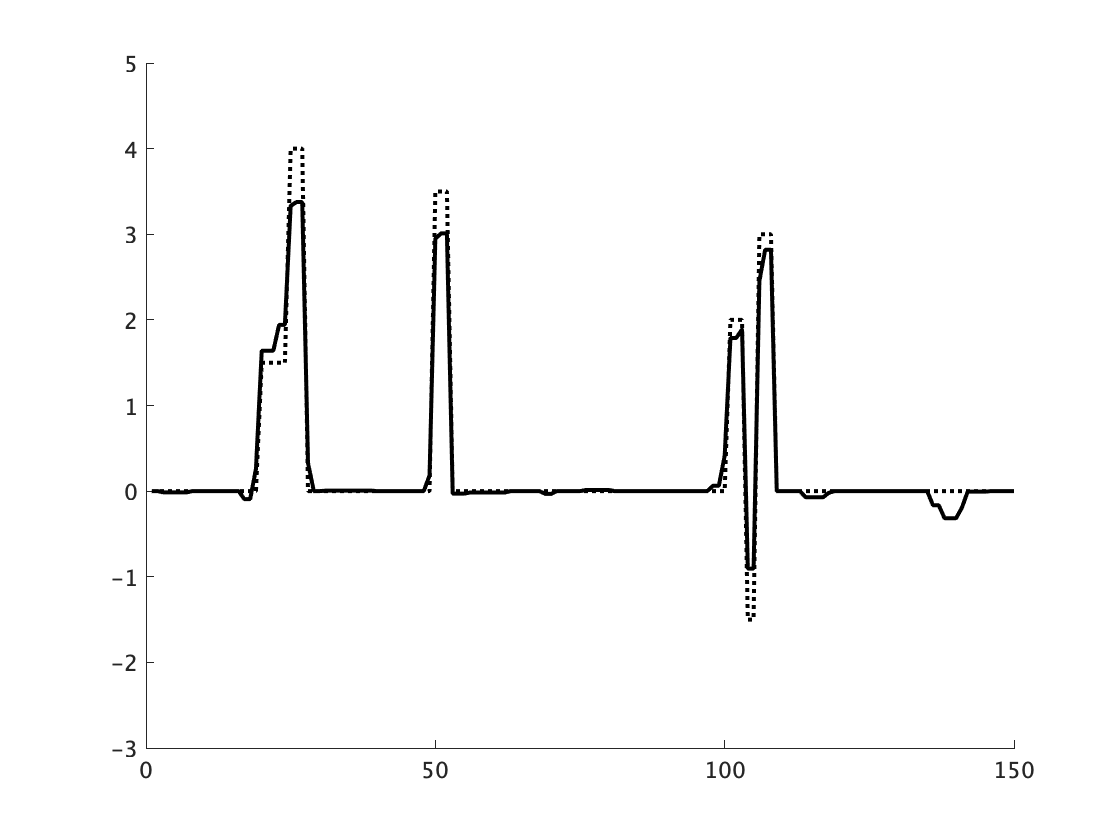}
    \subcaption{fused lasso (SE: $4.46$)}
  \end{minipage}
  \begin{minipage}[b]{0.46\linewidth}
    \centering
    \includegraphics[bb = 0 0 1120 840, scale=0.085]
    {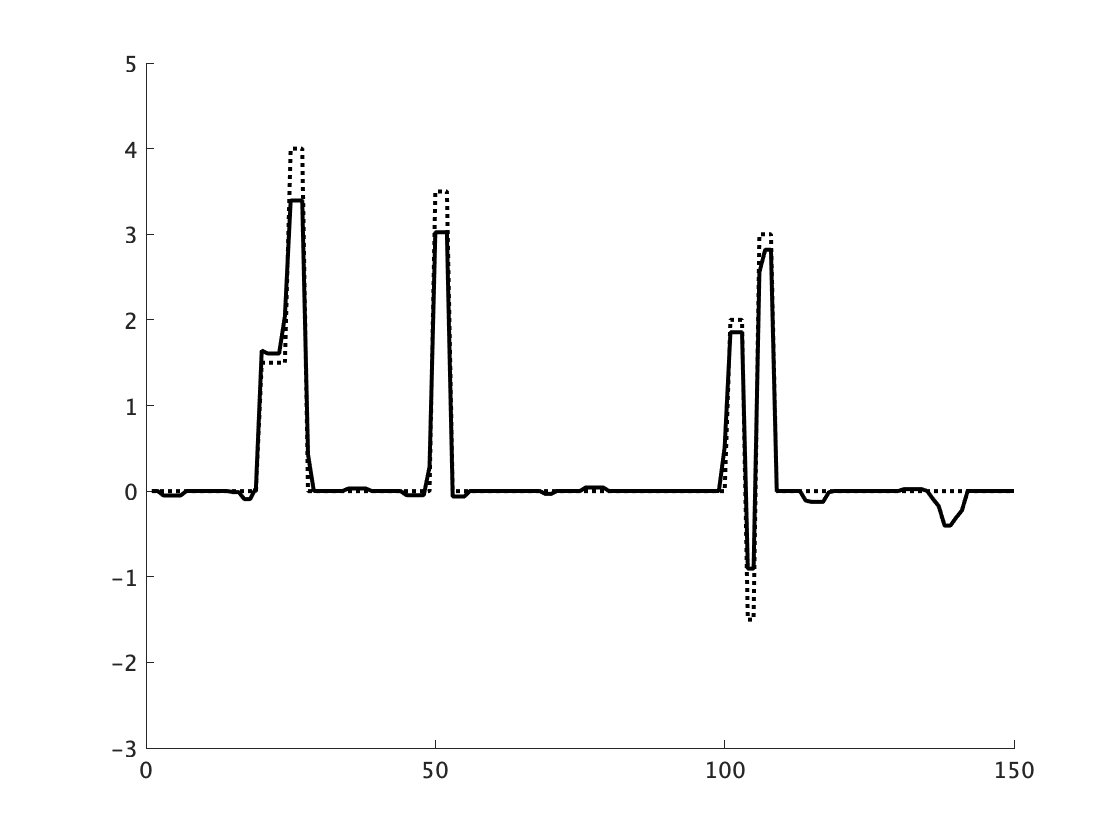}
    \subcaption{latent fused lasso (SE: $4.33$)}
  \end{minipage}\\
  \begin{minipage}[b]{0.46\linewidth}
    \centering
    \includegraphics[bb = 0 0 1120 840,scale=0.085]
    {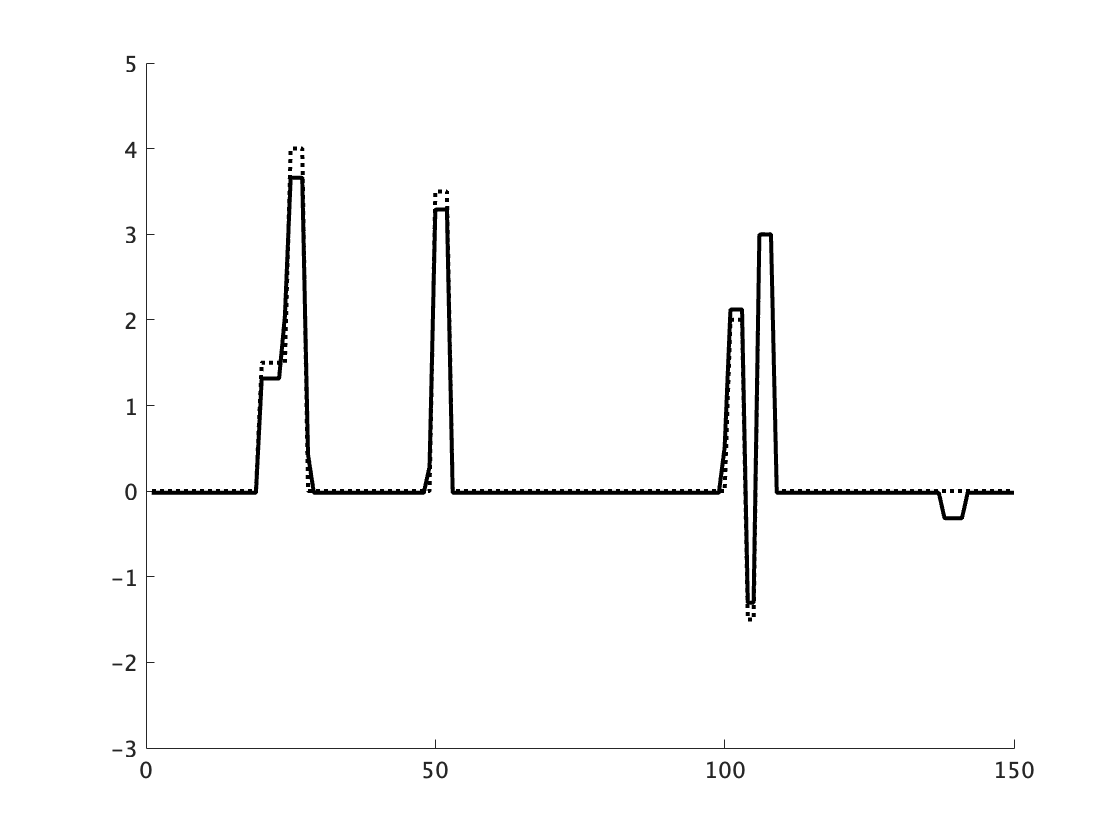}
    \subcaption{proposed model (i) (SE: $2.00$)}
  \end{minipage}
  \begin{minipage}[b]{0.46\linewidth}
    \centering
    \includegraphics[bb = 0 0 1120 840,scale=0.085]
    {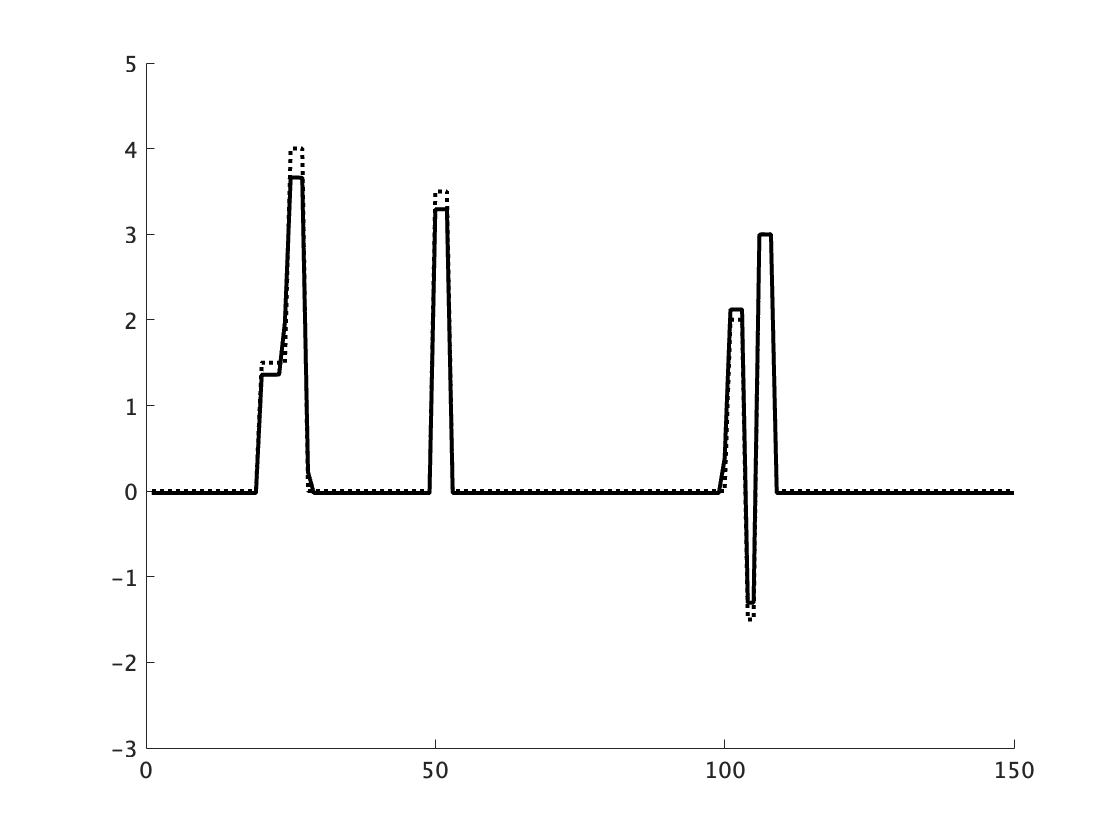}
    \subcaption{proposed model (ii) (SE: $ 1.15$)}
  \end{minipage}\\
  \begin{minipage}[b]{0.46\linewidth}
    \centering
    \includegraphics[bb = 0 0 1120 840,scale=0.085]
    {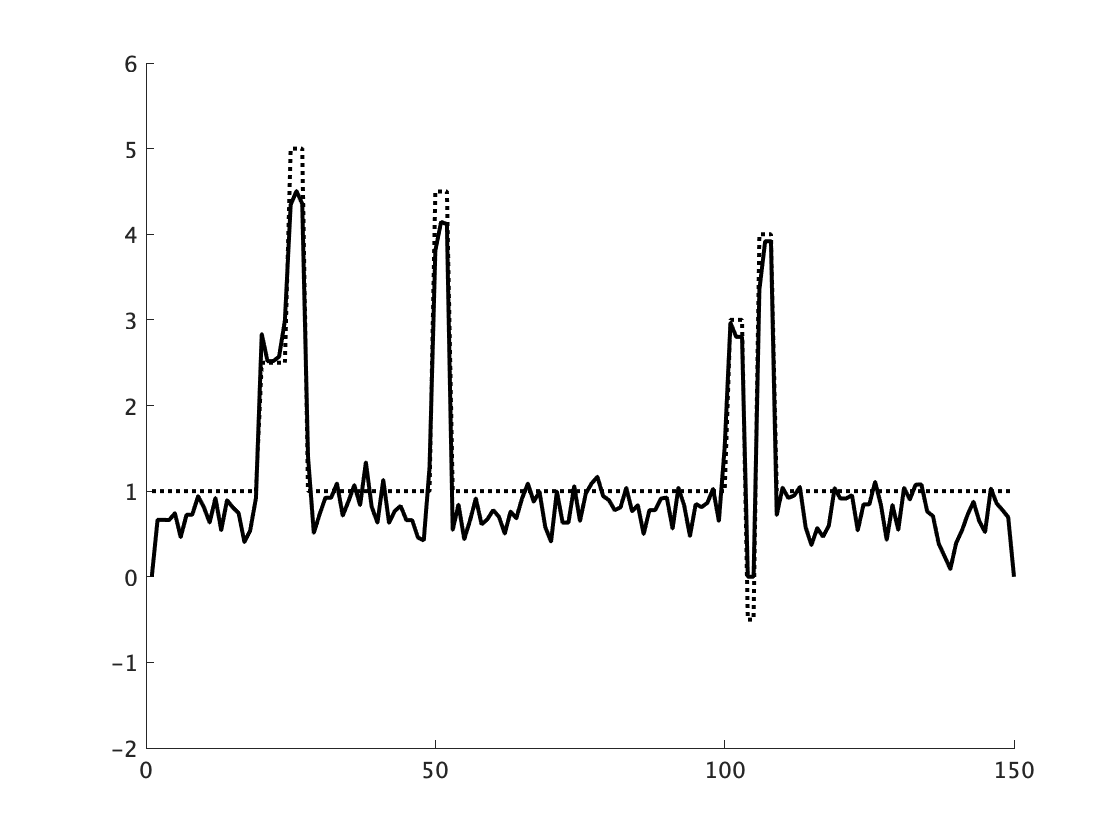}
    \subcaption{latent fused lasso}
  \end{minipage}
  \begin{minipage}[b]{0.46\linewidth}
    \centering
    \includegraphics[bb = 0 0 1120 840,scale=0.085]
    {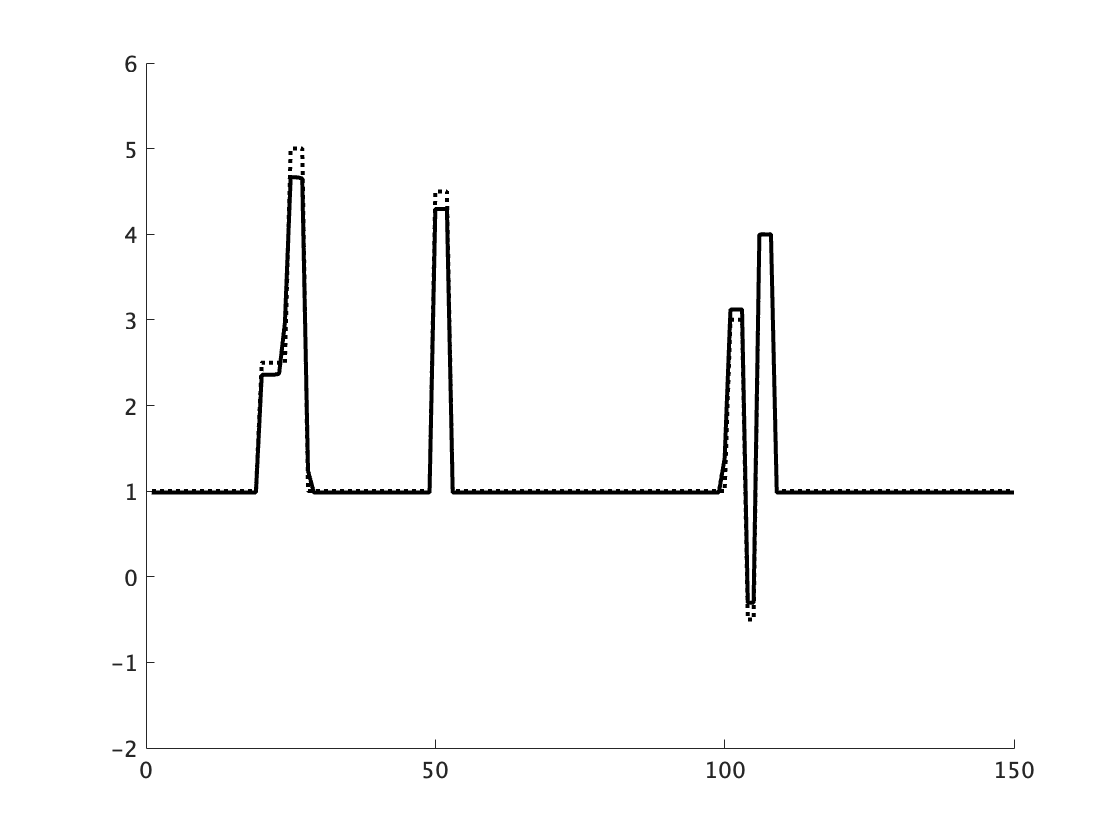}
    \subcaption{proposed model (ii)}
  \end{minipage}
  \caption{
  (a) original target signal $\xstar$,
  (b) observed noisy signal $y$,
  (c-f) denoised signals from $y$,
  (g-h) denoised signals from $y+\mathbf{1}_N$
  }
  \label{fig:denoise-example}
\end{figure}

We conducted numerical experiments in a scenario
of denoising of a sparse piecewise constant signal.
Fig.\ref{fig:denoise-example}(b) shows the observed signal $y=\xstar+ \varepsilon$, where $\xstar\in \setR^N\, (N=150)$ in Fig.\ref{fig:denoise-example}(a) is the target sparse piecewise constant signal to be estimated and $\varepsilon\in \setR^N$ is additive white Gaussian noise to achieve $10\log_{10}\frac{\norm{\xstar}^2}{\norm{\varepsilon}^2}= 10\mathrm{dB}$.
In the flame of Problem \ref{prob:unification-and-robustification} and \ref{prob:ne-fused-lasso}, we compared the denoised estimates by the fused lasso (see Remark \ref{remark:generalization-of-fls}(a)), by the latent fused lasso (see Remark \ref{remark:generalization-of-fls}(b)), by a nonconvex enhancement of the latent fused lasso realized by \eqref{eq:latent-fused-cLiGME} with $\mu_1=0$ (proposed model (i)), and by the general proposed model \eqref{eq:latent-fused-cLiGME} (proposed model (ii)).
For the proposed model (i) and (ii), we used $A=I_N$, $\cb^{(\mathcal{I})}=\cb^{(\setR)}$ and $\cd= \{\mathbf{0}_{\setR^{N-L}}\}$.
By Remark \ref{remark:why-nonconvex}(c), we obtained $B^{\lrangle{2}}$ for the proposed model (i) and $(B^{\lrangle{1}},B^{\lrangle{2}})$ for the proposed model (ii), where we used $(\theta_1, \theta_2) = (0.99, 0.99)$ for both models, $(\omega_1,\omega_2):=(0,1)$ for the proposed  model (i) and $(\omega_1,\omega_2):=(0.5,0.5)$ for the proposed model (ii).
For the latent fused lasso and the proposed models (i) and (ii), we set $L = 5$ which corresponds to the largest width of the non-baseline pulses in 
$\xstar$ (see also \cite[Remark]{feng2020latent}).
All denoised estimates except for those by the proposed model (ii) are obtained after $600$k iterations by
Algorithm \ref{alg:enforced}.
The denoised estimate by the proposed model (ii) is obtained after $1400$k iterations by Algorithm \ref{alg:enforced}.

Fig.\ref{fig:denoise-example}(c-e) show
the best denoised results $\widetilde{x}$ achieving the minimum squared error (SE): $\norm{\xstar- \widetilde{x}}^2$, after an exhaustive (but fairly efficient grid) search of regularization weights (i.e., $(\mu_1,\mu_2)$), by the fused lasso with $(\mu_1,\mu_2)=(0.1,0.3)$, by the latent fused lasso with $(\mu_1,\mu_2)=(0,0.4)$, and by the proposed model (i) with $(\mu_1,\mu_2)= (0,1)$. Fig.\ref{fig:denoise-example}(f) shows the best denoised estimate $\widetilde{x}$ by the proposed model (ii) with $(\mu_1,\mu_2)=(0.3,0.9)$ (for keeping $\frac{\mu_1}{\mu_2}=\frac{1}{3}$ in the above case of the fused lasso).
From Fig.\ref{fig:denoise-example}(d) and (e),
we see that the proposed nonconvex enhancement alleviates the underestimation effect caused by $\norm{\cdot}_1$ in the latent fused lasso.
From Fig.\ref{fig:denoise-example}(e) and (f), we see that \eqref{eq:latent-fused-cLiGME} with $\mu_2>0$ further improves the denoising performance.

To see the robustness against baseline mismatch, we also experimented after replacing the original target signal $\xstar$ in Fig.\ref{fig:denoise-example}(a) with $x^\star+ \mathbf{1}_N$ and the original observed noisy signal $y=\xstar+\varepsilon$ in Fig.\ref{fig:denoise-example}(b) with $y+\mathbf{1}_N$.
Fig.\ref{fig:denoise-example}(g) and (h) show the denoised estimates from $y+\mathbf{1}_N$, respectively, by the latent fused lasso 
with $\mathcal{I}=\{0\}$ (see Remark \ref{remark:generalization-of-fls}) according to zero baseline assumption (assumed implicitly or explicitly in \cite{feng2020latent}), and by the proposed model (ii) with $\mathcal{I}=\setR$ (both model use the best weights in the before experiment). 
These show the baseline robustness by $\cb^{(\mathcal{I})}$ with proper choice of $\mathcal{I}$.

\section{Conclusion}
We proposed the EA-type constrained algorithm for the cLiGME model.
As an application of the proposed algorithm,
we also proposed a bivariate nonconvex enhancement of a unified fused lasso models with effective constraints.
Numerical experiments show that the proposed model \eqref{eq:latent-fused-cLiGME} achieves superior estimation accuracy to conventional fused lasso models.

\bibliographystyle{IEEEbib}
\bibliography{references}

\begin{thebibliography}{10}

\bibitem{LiGME}
J.~Abe, M.~Yamagishi, and I.~Yamada,
\newblock ``Linearly involved generalized {Moreau} enhanced models and their
  proximal splitting algorithm under overall convexity condition,''
\newblock {\em Inverse Problems}, vol. 36, no. 3, pp. 035012, 2020.

\bibitem{selesnick2017}
I.~Selesnick,
\newblock ``{Sparse Regularization via Convex Analysis},''
\newblock {\em IEEE Transactions on Signal Processing}, vol. 65, no. 17, pp.
  4481--4494, 2017.

\bibitem{zhang2010}
C.-H. Zhang,
\newblock ``{Nearly unbiased variable selection under minimax concave
  penalty},''
\newblock {\em Annals of Statistics}, vol. 38, no. 2, pp. 894--942, 2010.

\bibitem{liu2022}
X.~Liu and E.C. Chi,
\newblock ``{Revisiting convexity-preserving signal recovery with the linearly
  involved GMC penalty},''
\newblock {\em Pattern Recognition Letters}, vol. 156, pp. 60--66, 2022.

\bibitem{Chen2023}
Y.~Chen, M.~Yamagishi, and I.~Yamada,
\newblock ``{A Unified Design of Generalized Moreau Enhancement Matrix for
  Sparsity Aware LiGME Models},''
\newblock {\em IEICE Transactions on Fundamentals of Electronics,
  Communications and Computer Sciences}, vol. E106.A, no. 8, pp. 1025--1036,
  2023.

\bibitem{yata2021}
W.~Yata, M.~Yamagishi, and I.~Yamada,
\newblock ``{A Constrained Linearly involved Generalized Moreau Enhanced Model
  and Its Proximal Splitting Algorithm},''
\newblock in {\em 2021 IEEE 31st International Workshop on Machine Learning for
  Signal Processing (MLSP)}, 2021, pp. 1--6.

\bibitem{cLiGME}
W.~Yata, M.~Yamagishi, and I.~Yamada,
\newblock ``{A constrained LiGME model and its proximal splitting algorithm
  under overall convexity condition},''
\newblock {\em Journal of Applied and Numerical Optimization}, vol. 4, no. 2,
  pp. 245--271, 2022.

\bibitem{combettes2021fixedpointstra}
P.~L. Combettes and J.-C. Pesquet,
\newblock ``Fixed point strategies in data science,''
\newblock {\em IEEE Transactions on Signal Processing}, vol. 69, pp.
  3878--3905, 2021.

\bibitem{shabili2021}
A.~H. Al-Shabili, Y.~Feng, and I.~Selesnick,
\newblock ``Sharpening sparse regularizers via smoothing,''
\newblock {\em IEEE Open Journal of Signal Processing}, vol. 2, pp. 396--409,
  2021.

\bibitem{zhang2023}
Y.~Zhang and I.~Yamada,
\newblock ``{A Unified Framework for Solving a General Class of Nonconvexly
  Regularized Convex Models},''
\newblock {\em IEEE Transactions on Signal Processing}, vol. 71, pp.
  3518--3533, 2023.

\bibitem{herman1980}
G.~T. Herman,
\newblock {\em Image Reconstruction from Projections},
\newblock Academic Press, 1980.

\bibitem{sezan1992}
M.~I. Sezan,
\newblock ``An overview of convex projections theory and its application to
  image recovery problems,''
\newblock {\em Ultramicroscopy}, vol. 48, no. 1, pp. 55--67, 1992.

\bibitem{byrne2004}
C.~Byrne,
\newblock ``A unified treatment of some iterative algorithms in signal
  processing and image reconstruction,''
\newblock {\em Inverse Problems}, vol. 20, no. 1, pp. 103, 2003.

\bibitem{censor2005}
Y.~Censor, T.~Elfving, N.~Kopf, and T.~Bortfeld,
\newblock ``The multiple-sets split feasibility problem and its applications
  for inverse problems,''
\newblock {\em Inverse Problems}, vol. 21, no. 6, pp. 2071--2084, 2005.

\bibitem{kitahara2021}
D.~Kitahara, R.~Kato, H.~Kuroda, and A.~Hirabayashi,
\newblock ``{Multi-Contrast CSMRI Using Common Edge Structures with LiGME
  Model},''
\newblock in {\em 2021 29th European Signal Processing Conference (EUSIPCO)},
  2021, pp. 2119--2123.

\bibitem{tibshirani2005}
R.~Tibshirani, M.~Saunders, S.~Rosset, J.~Zhu, and K.~Knight,
\newblock ``{Sparsity and Smoothness Via the Fused Lasso},''
\newblock {\em Journal of the Royal Statistical Society Series B: Statistical
  Methodology}, vol. 67, no. 1, pp. 91--108, 12 2005.

\bibitem{feng2020latent}
Y.~Feng and I.~Selesnick,
\newblock ``{Latent Fused Lasso},''
\newblock in {\em 2020 IEEE International Conference on Acoustics, Speech and
  Signal Processing (ICASSP)}, 2020, pp. 5969--5973.

\bibitem{angelosante2010}
D.~Angelosante, G.~B. Giannakis, and N.~D. Sidiropoulos,
\newblock ``{Estimating Multiple Frequency-Hopping Signal Parameters via Sparse
  Linear Regression},''
\newblock {\em IEEE Transactions on Signal Processing}, vol. 58, no. 10, pp.
  5044--5056, 2010.

\bibitem{lin2017}
K.~Lin, J.~L. Sharpnack, A.~Rinaldo, and R.~J. Tibshirani,
\newblock ``{A Sharp Error Analysis for the Fused Lasso, with Application to
  Approximate Changepoint Screening},''
\newblock in {\em Advances in Neural Information Processing Systems}, 2017,
  vol.~30.

\bibitem{kuroda2018}
H.~Kuroda, M.~Yamagishi, and I.~Yamada,
\newblock ``{Exploiting Sparsity in Tight-Dimensional Spaces for Piecewise
  Continuous Signal Recovery},''
\newblock {\em IEEE Transactions on Signal Processing}, vol. 66, no. 24, pp.
  6363--6376, 2018.

\bibitem{tibishirani2007}
R.~Tibshirani and P.~Wang,
\newblock ``{Spatial smoothing and hot spot detection for CGH data using the
  fused lasso},''
\newblock {\em Biostatistics}, vol. 9, no. 1, pp. 18--29, 05 2007.

\bibitem{mishra2023}
A.~Mishra, U.~K. Sahoo, and S.~Maiti,
\newblock ``Second-order fused lasso algorithm for radio tomographic imaging,''
\newblock {\em IEEE Communications Letters}, vol. 27, no. 7, pp. 1764--1768,
  2023.

\bibitem{tibshirani1996}
R.~Tibshirani,
\newblock ``{Regression Shrinkage and Selection via the Lasso},''
\newblock {\em Journal of the Royal Statistical Society. Series B
  (Methodological)}, vol. 58, no. 1, pp. 267--288, 1996.

\bibitem{CAaMOTiH}
H.~H. Bauschke and P.~L. Combettes,
\newblock {\em Convex Analysis and Monotone Operator Theory in Hilbert Spaces},
\newblock Springer, 2 edition, 2017.

\bibitem{yata2024}
W.~Yata and I.~Yamada,
\newblock ``{Imposing early and asymptotic constraints on LiGME with
  application to nonconvex enhancement of fused lasso models},''
\newblock in {\em 2024 IEEE International Conference on Acoustics, Speech and
  Signal Processing (ICASSP) (to appear)}, 2024.

\end{thebibliography}
\end{document}